\documentclass[12pt]{amsart}
\advance\hoffset by -0.5 truecm
\usepackage{amssymb}
\newtheorem{Theorem}{Theorem}[section]
\newtheorem{Lemma}[Theorem]{Lemma}
\newtheorem{Corollary}[Theorem]{Corollary}

\newtheorem{Remark}[Theorem]{Remark}

\def \dim{{\mbox {dim}}\,}
\def \ex{\mbox{\rm ex}}

\def\V{\mbox{Var}}

\def\Z{{\mathbb Z}}
\def\R\re
\def\V{\bf V}

\def \re{{\mathbb R}}
\def \Q{{\mathbb Q}}
\def \cp{{\mathbb CP}}

\def \0{\lambda_{0}}

\begin{document}
\title[Einstein manifolds and entropy]{Einstein manifolds of non-negative sectional curvature and entropy}

\author[G. P. Paternain]{Gabriel P. Paternain}\thanks{G.P. Paternain is on leave from Centro de Matem\'atica, Facultad de Ciencias, Igu\'a 4225, 11400 Montevideo, Uruguay.}
\address{CIMAT  \\
          A.P. 402, 36000 \\
          Guanajuato. Gto. \\
          M\'exico.}
\email {paternain@cimat.mx}

\author[J. Petean]{Jimmy Petean}
 \address{CIMAT  \\
          A.P. 402, 36000 \\
          Guanajuato. Gto. \\
          M\'exico.}
\email{jimmy@cimat.mx}

\thanks{J. Petean is supported by grant 28491-E of CONACYT}


\date{2000}


\begin{abstract}We show that if $(M^{n},g)$ is a closed simply connected Einstein manifold
of non-negative curvature then $-\log\,R\leq \frac{\pi\,\sqrt{n-1}\,(n-2)}{2}$
where $R$ is the radius of convergence of the series $\sum_{i\geq 2}\dim \,(\pi_{i}(M)\otimes\Q)t^{i}.$
If we suppose in addition that $M$ is formal then we show that:
$$\dim\, H_{*}(M,\Q)\leq \left[1+\exp\left(\frac{\pi\,\sqrt{n-1}\,(n-2)}{2}\right)\right]^{n}.$$
These results are achieved by combining the classical Morse theory of the loop space with a new
upper bound for the topological entropy of the geodesic flow of $g$ in terms of the curvature tensor.

\end{abstract}

\maketitle

\section{Introduction}Let $(M^{n},g)$ be a complete smooth connected
Riemannian manifold. One of our main concerns here
will be the study of Einstein metrics $g$ of non-negative sectional curvature.
Unless $(M,g)$ is flat, the Ricci curvature $r$ must be positive and we might as well
rescale the metric so that $r=g$. Therefore all the sectional curvatures satisfy $0\leq K(P)\leq 1$.
Myers' theorem asserts that $M$ is compact, has finite fundamental group and
its diameter is $\leq \sqrt{n-1}\,\pi$.
By passing to the universal covering we shall assume in the sequel that $M$ is
in fact closed and simply connected.

When $M$ is four-dimensional, strong topological restrictions to the existence of such metrics were first given
by N. Hitchin \cite{H}. He proved that the Euler characteristic $\chi$ and the signature
$\tau$ of $M$ must satisfy
\[\chi\geq \left(\frac{3}{2}\right)^{3/2}|\tau|.\]
This inequality was recently improved by M.J. Gursky and C. LeBrun \cite{GLe} who showed that in fact
one has $9\geq \chi>15/4|\tau|$. This result together with Freedman's classification
theorem \cite{Fr} tells us that there are at most 12 homeomorphism types (which can be explicitly
listed) of simply connected closed Einstein manifolds of non-negative sectional curvature.
If one assumes further that $M$ has positive definite intersection form, then
Gursky and LeBrun proved that $(M,g)$ must, up to isometry, be $\cp^{2}$ equipped
with a constant multiple of the Fubini-Study metric.
All these results rely heavily on the four-dimensional Gauss-Bonnet theorem.

Even though it is believed that simply connected closed Einstein manifolds
with non-negative sectional curvature are very rare, there are basically no known topological obstructions to
the existence of such metrics in dimensions $\geq 5$, except, of course, Gromov's celebrated result
in \cite{G1} which gives an obstruction to the existence of metrics just with non-negative sectional curvature.
Gromov's result asserts that if $M$ admits such a metric, then the sum of the Betti numbers $b_{i}$
must verify:
\[\dim\, H_{*}(M,\Q)=\sum_{i=0}^{n}b_{i}\leq C(n),\]
where 
\[C(n):=((n+1)J(n))^{100^{n}},\;\;\;\;\;J(n):=2^{M(n)},\;\;\;\;\;M(n)=8^{n}10^{n^{2}+4n}.\]
Clearly, the bound $C(n)$ is quite unrealistic, and Gromov himself conjectured
that $C(n)$ should be just $2^{n}$.

We now state our first main result. Let $R$ be the radius of convergence of the series
\[\sum_{i\geq 2}\dim \,(\pi_{i}(M)\otimes\Q)t^{i}.\]

\vspace{.5cm}

\noindent {\bf Theorem A.} {\it Let $(M^{n},g)$ be a smooth simply connected closed manifold
of positive Ricci curvature. Suppose that we normalize $g$ so that $r\geq g$
and let $k:=\max_{P} K(P)$. Then
\[-\log\,R\leq \frac{\pi\,\sqrt{n-1}}{2}\left((n-1)\sqrt{k}-\frac{1}{\sqrt{k}}\right).\]
In particular, if $(M^{n},g)$ is a simply connected closed Einstein manifold
of non-negative sectional curvature, then
\[-\log\,R\leq \frac{\pi\,\sqrt{n-1}\,(n-2)}{2}.\]}

\vspace{.5cm}
We remark that Theorem A also holds if we just assume that $M$ has finite fundamental group.
Also note that the second inequality in the theorem (as well as all the corollaries below) can also
be obtained just assuming the existence of a metric $g$ with $r\geq g$ and $K(P)\leq 1$
for all $P$.

A simply connected closed manifold $M$ is said to be {\it rationally elliptic}
if the rational homotopy $\pi_{*}(M)\otimes\Q$ is finite dimensional, i.e.
there exists a positive integer $i_{0}$ such that for all
$i\geq i_{0}$, $\pi_{i}(M)\otimes\Q=0$.
The manifold $M$ is said to be {\it rationally hyperbolic} if it is not
rationally elliptic (cf. \cite{FH1,FH2,GH1} and references therein).

It was proved by Y. F\'elix and S. Halperin \cite{FH} that if $M$ is rationally hyperbolic, then
 the integers $\rho_{i}=\sum_{j\leq i}\dim\,\pi_{j}(M)\otimes\Q$ grow
 exponentially in
 $i$ (i.e. there exist $C>1$ and a positive integer $k$ such that
 if $i>k$ then $\rho_{i}\geq C^{i}$).
Hence if $M$ is rationally
elliptic we have $R=+\infty$ and if $M$ is rationally hyperbolic $R<1$. It follows
that Theorem A is meaningful only if $M$ is rationally hyperbolic.

The ``generic" manifold is rationally hyperbolic; rational ellipticity
is a severely restrictive condition.
Examples of rationally elliptic manifolds are
simply connected homogeneous spaces \cite{Serre}, manifolds that admit a
codimension one compact action \cite{GH2}, Dupin hypersurfaces \cite{GH2}
and any known manifold that admits a Riemannian metric of non-negative
 sectional curvature.
A conjecture attributed to R. Bott states that any compact simply
connected manifold that admits a metric of non-negative sectional
curvature must be rationally elliptic (cf. \cite{GH1}).
It is known that if $M$ is rationally elliptic then
$\dim\, H_{*}(M,\Q)\leq 2^{n}$ \cite{GH1,H1} and hence Bott's conjecture implies Gromov's conjecture on the
optimal bound for the sum of the Betti numbers of $M$.

A manifold $M$ is said to be {\it formal}
if there exists
a morphism of differential graded algebras from
the minimal model of $M$ to $(H^{*}(M,\Q),0)$ that induces
an isomorphism in cohomology. The interest for this class
of spaces lies in the fact that for them all the rational
homotopy invariants of $M$ can be obtained from $H^{*}(M,\Q)$.
In \cite{DGMS}, it is shown that compact simply connected
K\"ahler manifolds are formal. Also,
any manifold with dimension $\leq 6$ is formal and, more generally,
any $p-1$-connected manifold of dimension $\leq 4p-2$ is formal \cite{NM}.

In Section 2 we shall explain how the results of Y. F\'elix and J.C. Thomas in \cite{FT}
combined with Theorem A yield the following two corollaries:

\vspace{.5cm}

\noindent {\bf Corollary 1.} {\it Let $(M^{n},g)$ be a smooth simply connected closed Einstein manifold
of non-negative sectional curvature. If $M$ is formal we have:
$$\mbox{\rm dim}\, H_{*}(M,\Q)\leq \left[1+\exp\left(\frac{\pi\,\sqrt{n-1}\,(n-2)}{2}\right)\right]^{n}.$$}

\vspace{.5cm}

\noindent {\bf Corollary 2.} {\it Let $(M^{n},g)$ be a smooth simply connected closed Einstein manifold
of non-negative sectional curvature. Suppose that $M$ is formal and $(p-1)$-connected.
Then
\[b_{p}\leq \frac{n}{p}\,\exp\left(\frac{p\,\pi\,\sqrt{n-1}(n-2)}{2}\right).\]}

\vspace{.5cm}

From Corollary 2, we obtain right away:

\vspace{.5cm}

\noindent {\bf Corollary 3.} {\it Let $(M^{5},g)$ be a smooth simply connected closed Einstein manifold
of non-negative sectional curvature.
Then
\[b_{2}\leq \frac{5}{2}\,\exp\left(6\,\pi\right).\]}

\vspace{.5cm}

Corollary 3 implies, for instance, that the connected sum of $k$ copies of $S^{3}\times S^{2}$
cannot admit an Einstein metric of non-negative sectional curvature for 
$k>\frac{5}{2}\,\exp\left(6\,\pi\right)$.

Certainly the bounds in Corollaries 1 and 3 are far better than Gromov's bound. 
On the other hand our bounds are not as good as those obtained by Gursky and LeBrun in \cite{GLe} for four manifolds.
The best lower bound for $-\log\, R$ in dimension four has been obtained by I. Babenko in \cite{Ba}. He proved that
\[1/R\geq \frac{b_{2}+\sqrt{b_{2}^{2}-4}}{2}.\]
Combining this bound with Theorem A we obtain:
\[\frac{b_{2}+\sqrt{b_{2}^{2}-4}}{2}\leq \exp(\pi\,\sqrt{3}).\]
This implies that $b_{2}\leq 230$ while Gursky and LeBrun's result implies that
$b_{2}\leq 7$.

The proof of Theorem A is based on ideas first introduced by M. Berger and R. Bott in \cite{BB}
and further developed by N. Grossman in \cite{Gro}. 
Given $x$ and $y$ in $M$ and $T>0$, define $n_{T}(x,y)$ as the number of 
geodesic arcs joining $x$ and $y$ with length $\leq T$.
For each $T>0$, the counting function $n_{T}(x,y)$ is finite and locally constant on an 
open full measure subset of $M\times M$, and integrable on $M \times M$. 
If $g$ has positive Ricci curvature, i.e., $r\geq \delta\,g$ with $\delta>0$, then 
we shall see in Section 2 that the Morse theory of the loop space yields

\begin{equation}
\frac{-\sqrt{\delta}\,\log\,R}{\pi\,\sqrt{n-1}} \leq \limsup_{T\rightarrow +\infty}\frac{1}{T}\log \int_{M}n_{T}(x,y)\,dy,
\label{key}
\end{equation}
for any point $x\in M$. G.P. Paternain explained in \cite{P1} how
 Yomdin's theorem \cite{Y} can be used to prove that
\begin{equation}
\limsup_{T\rightarrow +\infty}\frac{1}{T}\log \int_{M}n_{T}(x,y)\,dy
\leq h_{top}(g),  
\label{Pineq}
\end{equation}
where $h_{top}(g)$ denotes the topological entropy of the geodesic flow of $g$.
Combining (\ref{key}) and (\ref{Pineq}) yields:
\begin{equation}
\frac{-\sqrt{\delta}\,\log\,R}{\pi\,\sqrt{n-1}} \leq h_{top}(g).
\label{pb}
\end{equation}
Inequality (\ref{pb}) is also pointed out by I. Babenko in \cite[Theorem 4]{Ba}.
While Berger, Bott and Grossman estimated the average counting function by using
Rauch's comparision theorem we shall use dynamical ideas to estimate $h_{top}(g)$
from above. As we shall see below our bound is better than those
of Berger, Bott and Grossman and does not involve any lower bound on the sectional
curvature.
In Section 3 we shall prove the following result which will imply Theorem A
and has independent interest.

\vspace{.5cm}

\noindent {\bf Theorem B.} {\it Let $(M^{n},g)$ be a closed Riemannian manifold and let
$K_{max}$ be a positive upper bound for the sectional curvature. Then
\[h_{top}(g)\leq \frac{n-1}{2}\sqrt{K_{max}}-\frac{\min_{v\in SM}\,r(v)}{2\sqrt{K_{max}}},\]
where $SM$ is the unit sphere bundle of $M$ and $r(v)$ is the Ricci curvature
in the direction of $v\in SM$.}

\vspace{.5cm}

Let $k$ be a positive number such that $|K(P)|\leq k$ for all 2-planes $P$. Then, clearly
$r\geq -(n-1)k\,g$ and hence Theorem B gives
\[h_{top}(g)\leq \frac{n-1}{2}\sqrt{k}+\frac{n-1}{2}\sqrt{k}=(n-1)\sqrt{k}.\]
The latter inequality, which is certainly weaker than that of Theorem B, was first proved
by A. Manning in \cite{Ma2}.
The upper bound for $h_{top}(g)$ in Theorem B is probably the best that one can obtain
in terms of the $C^{0}$-norm of the curvature tensor. 
Note that the bound becomes sharp for all the space forms
(Manning's bound is not sharp for positively curved space forms).
Also observe that there is no hope to obtain an upper bound for $h_{top}(g)$ purely
in terms of bounds of the Ricci curvature. Indeed, the $K3$ surfaces admit Ricci flat
metrics, but $h_{top}(g)>0$ for any $C^{\infty}$ metric $g$ since a $K3$ surface
is rationally hyperbolic \cite{P,P1} (any simply connected four-manifold with second Betti
number strictly bigger than two is rationally hyperbolic).

We conclude this introduction by comparing the upper bound in Theorem B with the
bounds obtained by Berger and Bott in \cite{BB}
and Grossman in \cite{Gro}. Grossman's bound gives Berger and Bott's bound
if one assumes that the manifold has positive sectional curvature. The first
important difference between their bounds and ours is that ours does not involve
any lower bound on the sectional curvatures; this is certainly an advantage if one is
interested in upper bounds for $-\log\,R$ as in Theorem A.
 But even in the case of manifolds
of non-negative sectional curvature our bound is sharper as we now explain.
To make the comparison easier suppose that
our Riemannian metric $g$ has $0\leq K(P)\leq 1$ for all 2-planes $P$.
In this case our Theorem B gives 
$$h_{top}(g)\leq \frac{1}{2}\left(n-1-\min_{v\in SM}r(v)\right).$$
Under the same hypotheses on $g$ the inequality in Proposition 5.1 in \cite{Gro}
yields
\[ \limsup_{T\rightarrow +\infty}\frac{1}{T}\log \int_{M}n_{T}(x,y)\,dy
\leq \frac{2(n-1)}{\pi}\log\left(2+\frac{\pi}{2}\right).\]
which is certainly weaker than our bound even ignoring the term involving
the Ricci curvature since $\frac{2}{\pi}\log\left(2+\frac{\pi}{2}\right)$
is aproximately $0.8103$.


\section{Proof of Theorem A and Corollaries 1 and 2}

\subsection{Proof of Theorem A}

As in the introduction, given $x$ and $y$ in $M$ and $T>0$, let $n_{T}(x,y)$ be the number of 
geodesic arcs joining $x$ and $y$ with length $\leq T$.
For each $T>0$, the counting function $n_{T}(x,y)$ is finite and locally constant on an 
open full measure subset of $M\times M$, and integrable on $M \times M$ \cite{BB,P}.

Suppose now that the points $x$ and $y$ are not conjugate along any geodesic connecting them.
Let $\Omega(x,y)$ be the space of piecewise smooth paths $\alpha:[0,1]\to M$ with $\alpha(0)=x$
and $\alpha(1)=y$.
Given a non-negative integer $q$, let $i_{q}(x,y)$ be the number of geodesic arcs from
$x$ to $y$ with index $\leq q$. The Morse inequalities applied to the energy functional
on $\Omega(x,y)$ give right away \cite{Mil}:
\begin{equation}
\sum_{i=0}^{q}b_{i}(\Omega(x,y))\leq i_{q}(x,y).
\label{mt}
\end{equation}
It follows from the proof of Myers' theorem that if the Ricci curvarture $r$ satisfies
$r\geq \delta\,g$ for $\delta>0$, then any geodesic arc between $x$ and $y$ with length
at least $\pi\,\sqrt{\frac{n-1}{\delta}}\,q$ has index at least $q$. Hence for any non-negative
integer $q$ we have:
\begin{equation}
i_{q}(x,y)\leq n_{\pi\,\sqrt{\frac{n-1}{\delta}}\,q}(x,y).
\label{my}
\end{equation}
Combining (\ref{mt}) and (\ref{my}) we obtain
\[\sum_{i=0}^{q}b_{i}(\Omega(x,y))\leq n_{\pi\,\sqrt{\frac{n-1}{\delta}}\,q}(x,y),\]
and integrating this inequality with respect to $y$ gives:
\[\sum_{i=0}^{q}b_{i}(\Omega(x,y))\leq \frac{1}{\mbox{\rm Vol}(M)}\int_{M}n_{\pi\,\sqrt{\frac{n-1}{\delta}}\,q}(x,y)\,dy,\]
and hence
\begin{equation}
\limsup_{q\rightarrow +\infty}\frac{1}{q}\log \sum_{i=0}^{q}b_{i}(\Omega(x,y))\leq 
\pi\,\sqrt{\frac{n-1}{\delta}}\limsup_{T\rightarrow +\infty}\frac{1}{T}\log \int_{M}n_{T}(x,y)\,dy.
\label{bbg}
\end{equation}
This argument and the last inequality are taken from \cite{BB,Gro}.

Now let $R_{\Omega}$ be the radius of convergence of the Poincar\'e series:
\[\sum_{i\geq 0}b_{i}(\Omega(x,y),\Q)t^{i}.\]
Since this series always has infinitely many non-zero coefficients we clearly have:
\[-\log\,R_{\Omega}=\limsup_{q\rightarrow +\infty}\frac{1}{q}\log \sum_{i=0}^{q}b_{i}(\Omega(x,y),\Q).\]
On the other hand Babenko showed in \cite{Ba1} that if $M$ is rationally hyperbolic, then
$R=R_{\Omega}$ where $R$ is the radius of convergence of:
\[\sum_{i\geq 2}\dim \,(\pi_{i}(M)\otimes\Q)t^{i}.\]
Using (\ref{bbg}) we obtain inequality (\ref{key}) in the introduction, namely:
\begin{equation}
\frac{-\sqrt{\delta}\,\log\,R}{\pi\,\sqrt{n-1}} 
\leq \limsup_{T\rightarrow +\infty}\frac{1}{T}\log \int_{M}n_{T}(x,y)\,dy.
\label{key1}
\end{equation}
 In \cite{P1}, G.P. Paternain explained how
 Yomdin's theorem \cite{Y} can be used to prove that (see also \cite{P}):
\begin{equation}
\limsup_{T\rightarrow +\infty}\frac{1}{T}\log \int_{M}n_{T}(x,y)\,dy
\leq h_{top}(g),  
\label{Pineq1}
\end{equation}
where $h_{top}(g)$ denotes the topological entropy of the geodesic flow of $g$.
Combining (\ref{key1}) and (\ref{Pineq1}) yields:
\begin{equation}
\frac{-\sqrt{\delta}\,\log\,R}{\pi\,\sqrt{n-1}} \leq h_{top}(g).
\label{pb1}
\end{equation}
Theorem A is an immediate consequence of inequality (\ref{pb1}), Theorem B and the following observation:
if $M$ admits an Einstein metric $g$ with non-negative sectional curvature, then we can
rescale $g$ so that $r=g$ and $0\leq K(P)\leq 1$.

\qed

\subsection{Proof of Corollary 1}

We will use the following result of Y. F\'elix and J.C. Thomas \cite{FT}:

\begin{Theorem}Suppose that $M$ is formal and rationally
hyperbolic and let
$P_{M}$ be the Poincar\'e polynomial of $H_{*}(M,\Z)$.
Write $P_{M}(t)=\prod_{i=1}^{n}(t-z_{i})$.
Then
\[R\leq \min_{1\leq i\leq n}|z_{i}|.\]
\label{formal-5}
\end{Theorem}

\begin{Corollary} Let $M^{n}$ be a closed formal rationally hyperbolic manifold.
Then
$$\mbox{\rm dim}\,H_{*}(M,\Q)\leq \left(1+\frac{1}{R}\right)^{n}.$$
\label{suma}
\end{Corollary}

\begin{proof}By Poincar\'e duality if $z$ is a root of $P_{M}(t)$ then
$1/z$ is also a root of $P_{M}(t)$. Hence if we write
\begin{equation}
P_{M}(t)=\prod_{i=1}^{n}(t-z_{i})
\label{roots}
\end{equation}
it follows that
\begin{equation}
\min_{1\leq i\leq n}|z_{i}|=\frac{1}{\max_{1\leq i\leq n}|z_{i}|}.
\label{1/1}
\end{equation}
Let us set for brevity $A:=\max_{1\leq i\leq n}|z_{i}|$. From
(\ref{roots}) we get:
\[b_{i}\leq \left(^{n}_{i}\right)\,A^{i},\]
hence
\[\dim\, H_{*}(M,\Q)\leq \sum_{i=0}^{n}\left(^{n}_{i}\right)\, A^{i}=(1+A)^{n}.\]
By Theorem \ref{formal-5} and (\ref{1/1}) we have $A\leq 1/R$ which
yields
\[\dim\, H_{*}(M,\Q)\leq \left(1+\frac{1}{R}\right)^{n}\]
as desired.

\end{proof}

Corollary 1 follows from Theorem A and Corollary \ref{suma} if $M$ is rationally hyperbolic.
If $M$ is rationally elliptic we always have $\dim\, H_{*}(M,\Q)\leq 2^{n}$.
 
\qed

\begin{Remark}{\rm It seems interesting to observe that from inequality (\ref{pb1})
and Corollary \ref{suma} it follows that if $(M^{n},g)$ is a closed simply connected formal
manifold with $r\geq \delta\,g$, $\delta>0$, then
\[\dim\, H_{*}(M,\Q)\leq \left[1+\exp\left(\pi\,h_{top}(g)\,\sqrt{\frac{n-1}{\delta}}\right)\right]^{n}.\]
This suggests that the following smooth invariant of a manifold that admits a metric of 
positive Ricci curvature might be of interest:
\[h_{r}(M):=\inf_{\{g:\,r(g)\geq g\}}h_{top}(g).\]
From the last inequality, if $M$ is formal and admits a metric of positive Ricci curvature we
have
\[\dim\, H_{*}(M,\Q)\leq \left[1+\exp\left(\pi\,\sqrt{n-1}\,h_{r}(M)\right)\right]^{n}.\]}

\end{Remark}

\subsection{Proof of Corollary 2} In \cite{FT}, F\'elix and Thomas point
out the following corollary of Theorem \ref{formal-5}:

\begin{Corollary} If $M^{n}$ is a closed $(p-1)$-connected formal manifold then
\[R\leq \left(\frac{n}{p\,b_{p}}\right)^{1/p}.\]
\label{p-1}
\end{Corollary}

Corollary 2 follows right away from the Corollary \ref{p-1} and Theorem A.

\qed

\section{Proof of Theorem B}

Let $SM \subset TM$ be the unit sphere bundle of $M$. We will consider
the geodesic flow ${\Phi}_t$ of $g$ acting on  $SM$. 
Given any $\theta \in SM$ we will denote by $c_{\theta}$ the
geodesic with initial condition $\theta$. The expression
$\theta =(x,v)$ will mean that $v \in T_x M$.

Recall that the metric $g$ on $M$ induces a metric on $TM$ (the
Sasaki metric) and for any $\theta \in TM$ an orthogonal decomposition
of $T_{\theta}TM$ into horizontal and vertical parts: 
$T_{\theta}TM = H_{\theta} \oplus V_{\theta}$ (the vertical space is
of course the tangent space of the fiber). Recall also that the
differential of ${\Phi}_t$ has a nice expression in geometric terms;
given $\xi =(w_1 ,w_2 )\in T_{\theta}TM$, $d_{\theta}{\Phi}_t (\xi)
=(J_{\xi} (t), {\dot{J}}_{\xi} (t) )$, where $J_{\xi}$ is the Jacobi
field along $c_{\theta}$
with initial conditions $J_{\xi} (0)=w_1 $ and ${\dot{J}}_{\xi} (0)
=w_2 $.

If $\theta =(x,v)$ we will 
denote by $S({\theta})$ the orthogonal complement 
of $(v,0)$ in $T_{\theta}SM$. Equivalently, $S({\theta})$ is the orthogonal
complement of the subspace spanned by $(v,0)$ and $(0,v)$ in
$T_{\theta}M$. It is easy to see that the subspaces $S(\theta )$
are invariant through the differential of the geodesic flow.

\vspace{.2cm}

The
proof of Theorem B will be based on Przytycki's inequality \cite{Prrr} for the
topological entropy of the geodesic flow which we will
describe now.

Given two real vector spaces with inner product $V$ and $W$ of the
same dimension and a linear transformation $f:V\rightarrow W$ the
{\em expansion} of $f$, $\ex(f)$, is the supremum over all non-trivial subspaces
of $V$ of the absolute of the determinant of $f|_V$.
An important (and trivial) property of the expansion is that given two linear
maps $f$ and $g$ we have $\ex(f\,g)\leq \ex(f)\,\ex(g)$. This will be used
below.

Przytycki's inequality is the following:

$$h_{top}(g)\leq\liminf_{t\rightarrow \infty} \frac{1}{t} \log 
\int_{SM}\ex(d_{\theta}{\Phi}_t ) \ d\theta .$$

Here we consider $d_{\theta}{\Phi}_t$ as a map between $S(\theta )$
and $S({\dot{c}}_{\theta} (t))$.

\begin{Remark}{\rm We only need to use Przytycki's inequality
although it is actually true that for a $C^{\infty}$ Riemanniann metric
one has Ma\~n\'e's formula \cite{M2}:
\[h_{top}(g)=\lim_{t\rightarrow \infty} \frac{1}{t} \log 
\int_{SM}\ex(d_{\theta}{\Phi}_t ) \ d\theta= \lim_{T\rightarrow \infty} \frac{1}{T} \log 
\int_{M\times M}n_{T}(x,y)\,dxdy.\]
Moreover, O.S. Kozlovski \cite{Koz} has shown that Przytycki's 
inequality is an equality for arbitrary $C^{\infty}$ maps.}

\end{Remark}

We will find our upper bound for the topological entropy by studying
$\ex(d_{\theta}{\Phi}_t )$ for small values of $t$. 
Given any $\delta >0$ we
have:

$$h_{top}(g)\leq \liminf_{i\rightarrow \infty} \frac{1}{\delta i}
\log \int_{SM} \ex(d_{\theta}(({\Phi}_{\delta})^i )) \  d\theta. $$

Now suppose that for some $\delta >0$ 
and $\alpha >0$ we have that $\ex(d_{\theta}
{\Phi}_{\delta}) \leq \alpha$ for all $\theta$. Then it follows that

$$h_{top}(g) \leq \lim_{i\rightarrow \infty} \frac{1}{\delta i}
\log ({\alpha}^i \mbox{\rm Vol}(SM))=\frac{\log (\alpha )}{\delta}.$$

We will get an estimate of the form $\alpha 
=1+\beta \delta +O({\delta}^2)$. 
Then $\log (\alpha ) = \beta \delta + O({\delta}^2 )$ and since the 
last inequality holds for any $\delta >0$ we get

\begin{equation}
h_{top}(g) \leq \beta.
\end{equation}

\vspace{.2cm}

We will prove Theorem B by estimating $\beta$. Consider the polar
decomposition of $d_{\theta}{\Phi}_t$: 
$d_{\theta} {\Phi}_t =O_t (\theta ) \ L_t (\theta )$ where 
$O_t (\theta ):S(\theta )\rightarrow S({\dot{c}}_{\theta} (t))$ is
a linear isometry and $L_t (\theta )$ is a symmetric 
positive endomorphism of $S({\theta})$. 
Then $ex(d{\Phi}_t )=ex(L_t )$. Of course, $L_t$ can be given
explicitly: $L_t = ((d{\Phi}_t )^* \ (d{\Phi}_t ))^{1/2}$.

\vspace{.2cm}

Consider the map ${\mathcal R}:S({\theta}) \rightarrow S({\theta})$
which in the decomposition into horizontal and 
vertical parts  is given by ${\mathcal R}(w_1 ,w_2 )=
(w_2 ,-R(v,w_1 )v )$, where $R$ is the curvature tensor and
$\theta =(x,v)\in SM$.

\begin{Lemma} $((d{\Phi}_{\delta})^* \ (d{\Phi}_{\delta}) )^{1/2}=
Id + \frac{\delta}{2}({\mathcal R}+{\mathcal R}^* ) + O({\delta}^2 ).$

\end{Lemma}

\begin{proof}
Given any $\theta =(x,v)\in SM$ let $c_{\theta}$ be the geodesic with
initial condition $\theta$. Let $T_t$ be the parallel transport
along $c_{\theta}$ from $T_{c_{\theta} (0)}$ to $T_{c_{\theta} (t)}$.
Let $e_1 ,e_2 ,...,e_{n-1}$ be an orthonormal basis of 
$\{v\}^{\perp}\subset T_{x}M$ by eigenvectors
of the symmetric transformation $u\mapsto R(v,u)v $. Let 
$E_i$ be the parallel vector field along $c_{\theta}$ with initial
condition $e_i$. Given $\xi \in T_{\theta} SM$ we can write
$J_{\xi} (t)= \sum_{i=1}^{n-1}
a_i (t)E_i (t)$ for some smooth functions $a_i$.
Then,

$$J_{\xi} (\delta ) =
\sum_{i=1}^{n-1} (a_i (0)+\delta a_i ' (0))T_{\delta} E_i (0)
 \ \ +O \ ({\delta}^2 )$$

\noindent
and since ${\dot{J}}_{\xi} (t)=
\sum_{i=1}^{n-1} a_i ' (t)E_i (t)$, we have

$${\dot{J}}_{\xi}  (\delta )=
\sum_{i=1}^{n-1} (a_i ' (0) +\delta a_i '' (0))T_{\delta} E_i
(0) \ \ + \ O({\delta }^2 ).$$
Now, from the Jacobi equation, we get that 
$\sum_{i=1}^{n-1} a_i ''(0)e_i =-R(v,
J_{\xi} (0))v$. Therefore

\begin{align*}
(J_{\xi} (\delta ), &{\dot{J}}_{\xi} (\delta ))=\\
&\left( T_{\delta} \left(
\sum_{i=1}^{n-1} (a_i (0) +\delta a_i '(0) ) e_i \right),
 T_{\delta} \left(
\sum_{i=1}^{n-1} a_i '(0)e_i \ -\delta R(v, J_{\xi} (0))v \right)
\right)
+ O({\delta}^2 ),
\end{align*}

\noindent
i.e.

$$d{\Phi}_{\delta}
=(T_{\delta},T_{\delta}) (Id + \delta {\mathcal R})+O({\delta}^2 ).$$
Since $T_{\delta}$ is orthogonal we get

$$(d{\Phi}_{\delta})^* (d{\Phi}_{\delta}) = Id + \delta ({\mathcal R}+
{\mathcal R}^* )+O({\delta}^2 )$$

\noindent
and therefore,

$$((d{\Phi}_{\delta} )^* \ (d{\Phi}_{\delta}))^{1/2}=Id+
(\delta /2) ({\mathcal R}+{\mathcal R}^* )+O({\delta }^2 ).$$

\end{proof}

From the lemma we obtain that $\ex( d{\Phi}_{\delta} )= 
\ex \left( Id + \frac{\delta}{2} ({\mathcal R}+
{\mathcal R}^* ) \right) +O({\delta}^2 )$; and we are left to compute 
$\ex(Id +\frac{\delta}{2} ({\mathcal R}+
{\mathcal R}^* ) )$. This is a positive definite symmetric
endomorphism of $S(\theta )$ and therefore the expansion is the 
product of the eigenvalues which are greater than or equal to one, provided that there exists
at least one eigenvalue greater than or equal to one.

We can compute the
eigenvalues explicitly using the orthonormal basis
of $S(\theta )$ given by:
\[\{(e_{1},e_{1}),\dots,(e_{n-1},e_{n-1}),(e_{1},-e_{1}),\dots,(e_{n-1},-e_{n-1})\}.\]
 Suppose that $R(v,e_i )v ={\lambda}_i e_i$. Then it is easy to check that

$$ \left( Id + (\delta /2) ({\mathcal R}+
{\mathcal R}^* ) \right)
(e_i ,e_i )=\left( 1+ (\delta /2)(1-{\lambda}_i )
\right) (e_i ,e_i )$$

\noindent
and

$$\left( Id + (\delta /2) ({\mathcal R}+
{\mathcal R}^* ) \right) 
(e_i ,-e_i ) = \left( 1+(\delta /2)({\lambda}_i -1) \right)
 (e_i ,-e_i ).$$
If ${\lambda}_i \leq 1$ for all $i$ we get that 
$$\ex(Id + (\delta /2) ({\mathcal R}+
{\mathcal R}^* )) = \prod_{i=1}^{n-1} (1+(\delta /2) (1-{\lambda}_i ))
=1+ (\delta /2)(n-1-r(v)) + O({\delta}^2),$$

\noindent
where $r(v)$ denotes, as before, the Ricci curvature in the direction
of $v$.

Therefore if the sectional curvature of $g$ is bounded above by 1 we
get, from the previous discussion, 
that 
$$h_{top}(g)\leq \frac{1}{2} \left( n-1-\min_{v\in SM} r(v)\right).$$ 

For a general $g$, let $k$ be the maximum of all the sectional
curvatures. If $k>0$ the metric $g_k =k  g$ has sectional
curvature bounded above by one and from the previous observation
we get that 
$$h_{top}(g_k )\leq \frac{1}{2}
\left( n-1-\frac{\min_{v\in SM_g } r_g
  (v)}{k} \right). $$

But $h_{top}(g_k )= (1/ \sqrt{k}) \  h_{top}(g)$ and therefore we get

$$h_{top}(g)\leq \frac{1}{2} \left( \sqrt{k}(n-1)-\frac{\min_{v\in SM}
    r(v)}{\sqrt{k}}\right) .$$
 
This finishes the proof of Theorem B.

\qed

\begin{Remark}{\rm If the sectional curvature of $g$ is non-positive (i.e. if $k\leq 0$),
we get that for any positive number $\rho$,

$$h_{top}(g_{\rho})\leq \frac{1}{2}
\left( n-1-\frac{\min_{v\in SM_g } r_g (v)}{\rho} \right)$$
and therefore,

$$h_{top}(g)\leq  \frac{1}{2} \left( \sqrt{\rho }(n-1)-\frac{\min_{v\in SM}
    r(v)}{\sqrt{\rho }} \right) .$$

The sharpest inequality is obtained by taking $\rho =
\frac{-\min_{v\in SM}r(v)}{n-1}$ and hence we have:
\begin{equation}
h_{top}(g)\leq \sqrt{-(n-1)\min_{v\in SM}r(v)}.
\label{last}
\end{equation}

The Bishop comparison theorem implies right away that if $r\geq -(n-1)g$ then
$\lambda$, the exponential growth rate of volume of balls in the universal covering,
satisfies
\[\lambda\leq n-1.\]
But when the sectional curvature is non-positive, Manning proved \cite{Ma2} that $h_{top}(g)=\lambda$
and we recover inequality (\ref{last}).}

\end{Remark}

\end{document}